\documentclass[12pt]{amsart}
\usepackage{amsmath,amssymb,amscd,array, mathrsfs }

\usepackage{amsmath,amscd,amssymb,amsthm,array}
\usepackage{color}

\usepackage{amsmath,amssymb,mathrsfs,amsthm, tikz-cd,mathrsfs}
\usepackage{comment}
\def\url#1{\expandafter\s

\tring\csname #1\endcsname}

\def\mmat #1,#2,#3,#4,{\text{\small\arraycolsep=3pt $
\begin{pmatrix}#1&#2\\#3&#4\end{pmatrix}$}}

\usepackage{hyperref}

\usepackage{comments}
\newComments\SBe{Said}{blue}
\newComments\SBo{Sofiane}{blue}
\newComments\AM{Nacer}{blue}
\newComments\DL{DL}{red}
\newComments\QEh{QEh}{blue}

\def\mmat #1,#2,#3,#4,{\text{\small\arraycolsep=3pt $
\begin{pmatrix}#1&#2\\#3&#4\end{pmatrix}$}}

\usepackage{lscape}
\usepackage{tikz-cd}

\usepackage{DLdef1}
\usepackage{multicol}

\def\mmat #1,#2,#3,#4,{\text{\small\arraycolsep=3pt $
\begin{pmatrix}#1&#2\\#3&#4\end{pmatrix}$}}

\renewcommand {\ssbegin}[2][*]
 {\refstepcounter{subsection}%
\if#1*
\addcontentsline{toc}{subsection}{\thesubsection.\hskip 1pc #2}%
\else
\addcontentsline{toc}{subsection}{\thesubsection.\hskip 1pc #2. #1}%
\fi
 \def \secno {\gdef \secno {}{\ssecfont
\thesubsection.\hskip 2ex}%
 }%
 \begin{#2}}

\renewcommand {\sssbegin}[2][*]
 {\refstepcounter{subsubsection}
\if#1*
\addcontentsline{toc}{subsubsection}{\thesubsubsection.\hskip 1pc #2}%
\else
\addcontentsline{toc}{subsubsection}{\thesubsubsection.\hskip 1pc #2. #1}
\fi
 \def \secno {\gdef \secno {}{\ssecfont \thesubsubsection.\hskip 2ex}%
 }%
 \begin{#2}}

\renewcommand {\parbegin}[2][*]
 {\refstepcounter{paragraph}
\if#1*
\addcontentsline{toc}{paragraph}{\theparagraph.\hskip 1pc #2}%
\else
\addcontentsline{toc}{paragraph}{\theparagraph.\hskip 1pc #2. #1}
\fi
 \def \secno {\gdef \secno {}{\ssecfont \theparagraph.\hskip 2ex}%
 }%
 \begin{#2}}

\rmnameii{etr}{etr}
\rmnameii{evv}{ev}

\newcommand{\Z}{\mathbb{Z}}
\newcommand{\shh} {\underline{\llcorner\!\shortmid \!\lrcorner}}

\newcommand{\sh} {\llcorner\!\shortmid \!\lrcorner}

\title[On Zinbiel and Tortkara]{On Zinbiel and Tortkara superalgebras}

\author{Sofiane Bouarroudj}

\address{Division of Science and Mathematics, New York University Abu Dhabi, P.O. Box 129188, Abu Dhabi, United Arab Emirates.}
\email{sofiane.bouarroudj@nyu.edu}

\author{Farukh Mashurov }

\address {SDU University, Kaskelen, Kazakhstan.}
\address {Shenzhen International Center for Mathematics, Southern University of Science and Technology, Shenzhen, 518055, China.}
\email{f.mashurov@gmail.com }

\thanks{Corresponding author: Farukh Mashurov    was supported by the Science Committee of the
Ministry of Science and Higher Education of the Republic of Kazakhstan (Grant No. AP22683764)}
\thanks{Sofiane Bouarroudj was supported by the grant NYUAD-065. 
}

\keywords{Zinbiel superalgebra, Tortkara superalgebra, Lie element, algebraic classification.}
 \subjclass[2010] {17A30;  17A50; 17A70.}

\begin{document}

\maketitle

\begin{abstract}
In this paper, we study Zinbiel superalgebras and special Tortkara superalgebras, highlighting key differences between the super and the non-super setting. We present examples of Zinbiel superalgebras with Rota-Baxter operators and construct a basis for free Zinbiel superalgebras. Moreover, we establish a superalgebraic analogue of the Lie criterion for Zinbiel superalgebras. In contrast to the classical case, some homomorphic images of special Tortkara superalgebras on two generators are exceptional. Finally, we present a classification of all Tortkara superalgebras of dimensions 2 and 3.
\end{abstract}
\section{Introduction}

\subsection{Zinbiel superalgebras} The variety of Zinbiel algebras is a class of non-associative algebras related to Leibniz algebras. They were first introduced by Loday in \cite{L} as Koszul duals  of Leibniz algebras, and now they are intensively studied in the literature. Following Loday, an algebra with the identity 
$$ 
a (b c)-(a b) c-(b a) c=0, 
$$
is called a {\it (right)Zinbiel}  algebra. These algebras also appear in control theory (see \cite {Kaw} and references therin) and have been  studied by many authors.    Dzhumadil'daev and Tulenbaev proved that any finite-dimensional Zinbiel algebra over complex numbers is nilpotent. Later,  Nauryzbekova and Umirbaev demonstrated that any proper subvariety of the variety of Zinbiel algebras over a field of characteristic zero is nilpotent, see \cite{NU}. Recently, Towers \cite{Towers} proved that every finite-dimensional Zinbiel algebra over an arbitrary field is nilpotent, extending earlier results.

It is well-known that integrations are closely related to Zinbiel algebras.  
Recently, Chapoton considered the connection between Zinbiel algebras and multiple zeta values in \cite{Chap}, proposing new algebraic constructions and conjectures. Dotsenko and Shadrin established some connections between Zinbiel algebras and the operad of rational functions defined by Loday, see \cite{DS}. Lemay considered the relationship between the fundamental theorems of calculus and Zinbiel algebras. He showed that the category of Zinbiel  algebras is equivalent to the category of FTC-pair, see  \cite{Le} for more details.  
Zinbiel algebras with a derivation were considered in  \cite{KMS} and are related to derived algebras and Manin white product in theory of operad. In addition, Zinbiel algebras with a derivation can be used to construct  transposed Poisson algebras and generalized Poisson algebras \cite{DIM25}.

The symmetrization and skew-symmetrization of the Zinbiel multipication were studied in \cite{AK, D2, DIM19, K, L}. Loday in \cite{L} observed that any Zinbiel algebra under the anti-commutator product $x \circ y := x  y + y  x$ turns into a commutative and associative algebra. On the other hand, it was  proved in \cite{AK} that every nilpotent commutative algebra can be embedded into a Zinbiel algebra with respect to the anti-commutator operation.

\subsection{Tortkara superalgebras} Following A. Dzhumadil’daev  \cite{D2}, an anti-commutative algebra with the identity: 
\begin{equation}\label{id: Tortkara}
    (ab)(cd)+(ad)(cb)=J(a,b,c)d+J(a,d,c)b,
\end{equation}
where $J(x,y,z)=(xy)z+(yz)x+(zx)y$, is called a {\it Tortkara} algebra. Originally, a Tortkara algebra was defined in \cite{D2}, as minus-algebra $A^{(-)}$ of Zinbiel algebra. That is, any Zinbiel algebra  under the commutator product $[a, b] =
a b- ba$ is Tortkara algebra, see  \cite{D2}.

Both Zinbiel and Tortkara algebras also arise in the study of iterated-integrals and streamed informations in terms of half-shuffle and area operators, see \cite{DET1, DET2, DLPR, SDLPR}. Bremner \cite{B}  studied special identities in terms of the Tortkara triple product $[a, b, c] = [[a,b],c]$ on a free Zinbiel algebra. All identities of degree 5 and 7 were obtained, which are not derived from  identities of lower degree. Bremmer also showed that, although Tortkara algebras are defined by a cubic binary
operad (with no Koszul dual), the corresponding triple systems are defined by
a quadratic ternary operad (with a Koszul dual). Algebraic and geometric classifications of Tortkara algebras in low dimensions  have been obtained in \cite{GKK18, GKK19}.  Kolesnikov in \cite{K} showed that the class of all special Tortkara algebras does not form a variety. Dzhumadil’daev and his coauthors in \cite{DIM19} established a criterion for elements of free Zinbiel algebras to be Lie. This criterion is applied to the study of speciality problems in Tortkara algebras, as well as to the construction of a basis for free special Tortkara algebras. In addition, they obtained an analogue of Cohn’s theorem on speciality of homomorphic
images of special Jordan algebras in two generators, see \cite{Cohn}. That is, any homomorphic image of a free special Tortkara algebra on two generators is special, see \cite{DIM19}. Cohn's theorem for superalgebras is also discussed in \cite{G}.

 Recently, some results in the theory of Zinbiel algebras have been  generalized to Zinbiel superalgebras \cite{BKM, COKN}. For instance, it has been proved that all Zinbiel superalgebras over an arbitrary field are nilpotent, analogously to Zinbiel algebras. Zinbiel superalgebras with low dimensions were classified in \cite{COKN}. 

 \subsection{Main results.}
The main objective of this paper is to continue the study of Zinbiel and special Tortkara in the super setting,  pointing out some differences between the non-super and super case. First, we introduce a new class of superalgebras called {\it Tortkara} superalgebras as these superalgebras verify the super Tortkara identity  Eq. (\ref{id: sup-tortkara}) which is a superization of the Tortkara identity Eq. (\ref{id: Tortkara}).  In Section \ref{examplesZ}, we provide some examples of Zinbiel superalgebras equipped with a Rota-Baxter operator, and construct a basis for free Zinbiel superalgebras. Furthermore, we establish an analogue of the Lie criterion for Zinbiel superalgebras, see Theorem \ref{Criterion}. Notably, we demonstrate that, unlike in the non-super setting, certain homomorphic images of special Tortkara superalgebras on two generators are exceptional, see Theorem \ref{hom image}. Finally, Section \ref{class} is devoted to the classification of Tortkara superalgebras of dimensions 2 and 3.

All our (super)algebras are considered over a field $\mathbb{K}$ of characteristic 0.

\section{Main Definitions} 

We denote the group of integers modulo $2$ by $\Z_{2}$. Let $V=V_\ev \oplus V_\od$ be a $\Z_{2}$-graded space. The degree of a homogeneous element $v\in V_{\bar{i}}$ is denoted by $|v|:=\bar{i}$. The element $v$ is called \textit{even} if $v\in V_\ev$ and \textit{odd} if $v\in V_\od$.  We denote by $\Pi$ the \textit{change of parity functor} $\Pi: V\mapsto \Pi (V)$, where $\Pi(V)$ is another copy of $V$ such that $ \Pi(V)_\ev:=V_\od;~~\Pi(V)_\od:=V_\ev$.  Elements of $\Pi(V)$ shall be denoted by $\Pi(v),~\forall v\in V$. We also shall identify $V$ and $\Pi(\Pi(V))$ in a natural way.\\

A superalgebra over $\mathbb K$ is a $\mathbb K$-module $A$ with a direct sum decomposition $A=A_\ev\oplus A_\od$ together with a bilinear multiplication $\cdot:A\times A \rightarrow A$. We usually write $ab$ instead of $a\cdot b$ unless we introduce a new operation. 


Recall that the super commutator and super anti-commutator defined as follows
(for all homogeneous $a,b$ in $A$)
\[
\begin{array}{lcl}
\{a , b \} & := & ab+(-1)^{|a||b|}ba, \\[1mm]
[a,b] & := & ab- (-1)^{|a||b|}ba.
\end{array}
\]

In the super setting, applying the Sign Rule sometimes requires dexterity. There are two versions of skew- or anti-commutativity that are synonyms only in the non-super case; for two homogeneous elements $a$ and $b$ of a superalgebra we call the following conditions
\[
\begin{array}{lcll}
ba & = & (-1)^{|b||a|}ab & \text{super commutativity,} \\[2mm]
ba & = & -(-1)^{|b||a|}ab & \text{super anti-commutativity, }\\[2mm]
ba & = & (-1)^{(|b|+1)(|a|+1)} ab & \text{super skew-commutativity, }\\[2mm]
ba & = &  -(-1)^{(|b|+1)(|a|+1)}ab & \text{super antiskew-commutativity.}
\end{array}
\]
For instance, the bracket in any Lie superalgebra is super anti-commutative while the Buttin bracket (see, e.g. \cite{GPS}) 
is, however, super antiskew-commutative relative to the natural parity of generating functions.


\subsection{Zinbiel and Tortkara superalgebras} \label{Grassmann}
We will build Zinbiel superalgebra and Tortkara superalgebra using the Grassman envelope method; see \cite{S, ZS}. 

Let $A = A_{\bar 0} \oplus A_{\bar 1}$ be a superalgebra (a $\Z_2$-graded algebra). Let $G = G_{\bar 0} \oplus G_{\bar 1}$ be the Grassmann algebra over a countable set of generators $\{ \xi_i | \, \xi_i\xi_j = -\xi_j\xi_i \}$ with the standard $\Z_2$-gradation in which $G_{\bar 0}$ (respectively $G_{\bar 1})$ is spanned by the words of even (respectively odd) length in the generators $\{\xi_i\}$. 

The {\it Grassmann envelope} of the superalgebra $A$ is the algebra
$$
G(A)=G_{\bar 0}\otimes A_{\bar 0}+G_{\bar 1}\otimes A_{\bar 1}.
$$

If $\mathcal{P}$ is some homogeneous variety of algebras and $G(A)$ belongs to the variety $\mathcal{P}$, then $A$ is called a $\mathcal{P}$-superalgebra. We will apply this method to superize the identity \eqref{id: Tortkara}.

A superalgebra $A=A_{\bar 0}\oplus A_{\bar 1}$ is called a Zinbiel superalgebra if it satisfies the following identity: (for all homogeneous $a,b,c$ in $A$) 
\begin{equation} \label{id: superZinb}
a (b c)-(a b) c-(-1)^{|a||b|}  ( ba) c=0.
\end{equation}

A superalgebra $T=T_{\bar 0}\oplus T_{\bar 1}$ is said to be a Tortkara superalgebra if and only if $G(T)$ is a Tortkara algebra. Direct computations reveal that this is equivalent to satisfying the following identities: (for every homogeneous elements $a,b,c,d$)
\begin{eqnarray}\label{id: sup-anti-com}
    ab & = & -(-1)^{|a||b|}ba,\\
\label{id: sup-tortkara}
    (ab)(cd)-(-1)^{|d|(|b|+|c|)}(ad)(bc) & = & J_s(a,b,c)d+(-1)^{|a||b|}b J_s(a,c,d),
\end{eqnarray}
where $J_s(x,y,z):=(xy)z-x(yz)-(-1)^{|y||z|}(xz)y$ (the super-Jacobian).

The expanded form of the identity \eqref{id: sup-tortkara} is shown below:
\[
\begin{array}{lcl}
(ab)(cd)-(-1)^{|b||d|+|c||d|}(ad)(bc) & = &((ab)c)d-(a(bc))d-(-1)^{|b||c|}((ac)b)d+(-1)^{|a||b|}b((ac)d)\\[2mm]
& & -(-1)^{|a||b|}b(a(cd))-(-1)^{|a||b|+|c||d|}b((ad)c).
\end{array}
\]

We will show that any Zinbiel superalgebra turns into a Tortkara superalgebra under the super-commutator operation.

A Tortkara superalgebra $T$ is called \textit{special} if there exists a Zinbiel superalgebra $A$ such that $T$ is a super-subalgebra of $A^{(-)}=(A,[\cdot,\cdot])$; otherwise, it is called \textit{exceptional}.


\section{Examples of Zinbiel superalgebras}\label{examplesZ}

The classical Rota–Baxter operators have been widely studied in Algebra, see, e.g., \cite{EG}. Recently, an odd version of this operator has been introduced in \cite{BGZ23}.  The generalization {of Rota-Baxter operators} in the context of superalgebra is given as follows:

Let $(A, \cdot)$ be a superalgebra. A homogeneous linear map $R: A \to A$ is called a Rota–Baxter operator of parity $|R|$ if it satisfies the identity  
$$R(x)R(y)=R((-1)^{|R|(|x|+|R|)}R(x)y+xR(y)),$$for all homogenous elements $x, y \in A$.

\ssbegin{Proposition}
     Let $(A,\cdot)$ be a supercommutative associative superalgebra and let $R$ be a Rota-Baxter operator. Let us define a product $\circ$ on $A$ as follows $$a\circ b := R(a)b.$$
     Then
     
     \begin{itemize}
         \item[(i)] If $R$ is an even Rota-Baxter operator, then $(A,\circ)$ is a Zinbiel superalgebra.
         \item[(ii)] If $R$ is an odd Rota-Baxter operator, then $(A,\circ)$ is a superalgebra with identity
         \begin{equation}\label{id: odd Zinbiel} a\circ (b\circ c)=(-1)^{|a|+1} (a\circ b)\circ c + (-1)^{|a|(|b|+1)} (b\circ a)\circ c.\end{equation}
     \end{itemize}
\end{Proposition}
\begin{proof}
    For every homogeneous  $a,b,c\in A$ we have 
    $$\begin{array}{ll}
       a\circ (b\circ c)  & =  R(a)R(b)c
       \\ & =R((-1)^{|R|(|a|+|R|)}R(a)b+a R(b)) c \\
         & =  R((-1)^{|R|(|a|+|R|)}R(a)b+(-1)^{|a|(|b|+|R|)}R(b)a) c\\
         & =(-1)^{|R|(|a|+|R|)} (a\circ b)\circ c + (-1)^{|a|(|b|+|R|)} (b\circ a)\circ c.
    \end{array}$$
    This completes the proof.
\end{proof}

As it turns out, the superalgebra satisfying \eqref{id: odd Zinbiel} is not new. Using the change of parity functor, we will demonstrate how it can be constructed from a Zinbiel superalgebra.

Let $Z=Z_{\bar 0}\oplus Z_{\bar 1} $ be a Zinbiel superalgebra with a product $\circ$. Consider the graded space $A:=\Pi(Z)$, where $\Pi$ is the changing parity functor. We define a new product $\star$ on $A$ as follows:
\[
a\star b:= (-1)^{|a|}\; \Pi( \Pi(a)\circ \Pi(b)). 
\]
Let us show that the product $\star$ satisfies \eqref{id: odd Zinbiel}. Indeed, for all homogeneous $a,b,c \in A$:
\[
\begin{array}{lcl}
a\star (b \star c) & = & (-1)^{|a|+|b|} \, \Pi (\Pi(a)\circ  (\Pi(b)\circ \Pi(c)))\\[2mm] 
&= & (-1)^{|a|+|b|} \, \Pi( ( \Pi(a)\circ \Pi(b)+ (-1)^{(|a|+1)(|b|+1)}\Pi(b)\circ \Pi(a)) \circ \Pi(c))
\\[2mm] 
&= & (-1)^{|b|} \, \Pi( (\Pi(a \star  b))\circ \Pi(c))+ (-1)^{(|a|+1)(|b|+1)+|a|}\Pi((\Pi(b\star a)) \circ \Pi(c))
\\[2mm] 
&= & (-1)^{|a|+1} \, (a \star  b)\star c+ (-1)^{|a|(|b|+1)}  (b\star a) \star c.
\end{array}
\]
Similarly, any superalgebra satisfying \eqref{id: odd Zinbiel} yields a Zinbiel superalgebra after changing  the parity and product; we omit the details.

Dzhumadil'daev constructed an interesting series of examples of Zinbiel and Tortkara algebras on a commutative associative algebra with integration operation \cite{D3}. We generalize these examples to the case of superalgebras. We begin by proving the following proposition.

\ssbegin{Proposition}\label{prop: Zin new is Zin} Let $(A,\circ)$ be a Zinbiel superalgebra and let $R$ be an even Rota-Baxter operator on $A$. Define the product $\circ_{1}$ on $A$ as follows $$a\circ_{1} b= R(a)\circ b+ a \circ R(b).$$ Then,  $(A,\circ_{1} )$ is a Zinbiel superalgebra.
\end{Proposition}
\begin{proof} Let $a,b,c\in A$. By definition of Rota-Baxter operator, we have
\begin{align*} a\circ_1(b \circ_1 c)=& R(a)\circ (R(b)\circ c +b\circ R(c))+a\circ R(R(b)\circ c+b\circ R(c))\\
=&  R(a)\circ (R(b)\circ c +b\circ R(c))+ a\circ (R(b)\circ R(c))\\
=&  (R(a)\circ R(b))\circ c +(-1)^{|a||b|}(R(b)\circ R(a))\circ c
 +(R(a)\circ b)\circ R(c)\\
 & +(-1)^{|a||b}(b\circ R(a))\circ R(c) + (a\circ R(b))\circ R(c)+(-1)^{|a||b}(R(b)\circ a)\circ R(c).
\end{align*}
Moreover,
\begin{align*} (a\circ_1 b) \circ_1 c=& R(R(a)\circ b+ a \circ R(b))\circ c+(R(a)\circ b+ a \circ R(b))\circ R(c)\\
=& (R(a) \circ R(b))\circ c+(R(a)\circ b)\circ R(c)+ (a \circ R(b))\circ R(c).
\end{align*}
Therefore, the identity \eqref{id: superZinb} holds.
\end{proof}

It is easy to see that $R$ is an even Rota-Baxter operator on the Zinbiel superalgebra $(A,\circ_1).$ Indeed,

$$\begin{aligned} R(a)\circ_1 R(b)=&R(R(a))\circ R(b)+R(a)\circ R(R(b))\\
=&R(R(R(a))\circ b +R(a)\circ R(b))+R(R(a)\circ R(b)+a\circ R(R(b))\\
=&R(R(a)\circ_1 b+a\circ_1 R(b)).
\end{aligned}$$

Then the above proposition can be used to define a series of Zinbiel superalgebras with an even Rota-Baxter operator $R$ on $A$.

\ssbegin{Corollary}\label{cor: series of Zinbiel} Let $(A,\circ)$ be a Zinbiel superalgebra and let $R$ be an even Rota-Baxter operator on $A$. Define inductively the product $\circ_{i+1}$ on $A$ as follows $$a\circ_{i+1} b= R(a)\circ_i b+ a \circ_i R(b),$$ where $\circ_0=\circ$ and $i\in \{0,1,\ldots,n\}.$ Then $(A,\circ_{i+1})$ is a Zinbiel superalgebra for every $i\in \{0,1,\ldots,n\}.$
\end{Corollary}

Furthermore, Dzhumadil'daev's examples in \cite{D3} can be interpreted as follows.

\ssbegin{Corollary} Let $(A,\cdot)$ be a supercommutative associative superalgebra and $R$ be a Rota-Baxter operator. Define the product $\circ$ on $A$ as follows $$a\circ_n b=\sum_{i=0}^n\binom{n}{i}R^i(a) R^{n-i+1}(b).$$
Then $(A,\circ_{i+1})$ is a Zinbiel superalgebra.    
\end{Corollary}
\begin{proof}  It is easy to see that $a\circ_n b=R(a)\circ_{n-1} b+ a \circ_{n-1} R(b)$ in the supercommutative associative superalgebra $(A,\cdot)$.   Corollary \ref{cor: series of Zinbiel} completes the proof.
\end{proof}

\section{A free Zinbiel superalgebra} \label{a free Zinbiel}

Let  $V = V_0\oplus V_1$ be a $\Z_2$-graded vector space over a field $ \mathbb{K}$ such that $V_{\overline{0}}$ is the space   of even elements, with parity $|x| = \bar 0$  and $V_{\overline{1}}$  is the space of odd elements, with parity  $|x| = {\bar 1}.$

The tensor superalgebra $ T(V) $ on the graded vector space $ V $ is defined as:
\[
T(V) = \bigoplus_{n \geq 0} V^{\otimes n},
\]
where $ V^{\otimes n} $ is the $ n $-fold tensor product of $ V, $ with $ V^{\otimes 0} = \mathbb{K}.$ Denote by $ v_1  v_2  \cdots v_n $ the generators $ v_1 \otimes v_2 \otimes \cdots \otimes v_n \in  V^{\otimes n} .$ Elements of $ T(V) $ are finite sums of tensors of the form $ v_1  v_2  \cdots  v_n $, where $ v_i \in V $. Each element inherits a grading from the $\Z_2$-grading of $ V $, with the total degree of $ v_1  v_2  \cdots v_n $ defined as:
\[
|v_1  v_2 \cdots  v_n| = |v_1| + |v_2| + \cdots + |v_n| \mod 2.
\]

Let $Sh_{m,n}$ be a set of shuffle permutations, i.e.,
$$Sh_{m,n}=\{\sigma\in S_{m+n} | \sigma(1)<\cdots<\sigma(m), \sigma(m+1)<\cdots<\sigma(m+n)\}.$$
For any positive integers  $i_1,\ldots,i_m,j_1,\ldots,j_n$ denote by 
$Sh(i_1\ldots i_m;j_1\ldots j_n)$ the set of  sequences $\alpha=\alpha_1\ldots\alpha_{n+m}$ constructed by 
shuffle permutations $\sigma\in Sh_{m,n}$  by changing $\alpha_{\sigma(l)}$ to $i_l$ if $l\le m$ and to $j_{l-m}$ if $m<l\le m+n.$

Let $\sigma\in Sh(i_1\ldots i_m;j_1\ldots j_n).$ We define  \begin{equation}\label{func: phi}K(\sigma)=\{(i,j)| \; i<j, \sigma_i>\sigma_j\}\quad \text{ and } \quad \phi(\sigma)=\sum_{(i,j)\in K(\sigma)} |v_{\sigma_i}||v_{\sigma_j}|.\end{equation}

Let us define the Grassmann envelope of the tensor algebra $T(V)$ by 
\[
G(T(V))=G_{\bar 0} \otimes T(V)_{\bar 0}\oplus G_{\bar 1} \otimes T(V)_{\bar 1},
\] 
where $G = G_{\bar 0} \oplus G_{\bar 1}$ is the Grassmann algebra over a countable set of generators $\{ \xi_i | \, \xi_i\xi_j = -\xi_j\xi_i \}$, see Section \ref{Grassmann}. Let $v_1, v_2, \dots$ be a homogeneous basis of $V$. Consider a basis $\hat{v}_1,\hat{v}_2,\ldots $ with the property that  
\[
\hat{v}_j = \left \{ 
\begin{array}{lcl}
1 \otimes v_j, & \text{if} &   |v_j|=\bar 0,\\
\xi_j \otimes v_j, & \text{if} &  |v_j|=\bar 1.
\end{array}
\right .
\]

The (ordinary) \emph{shuffle product} in the algebra  (see, e.g., \cite{DIM19, L}) is given by:
\begin{equation}\label{eq:ordinary_shuffle}
(\hat{v}_{i_1} \cdots \hat{v}_{i_p}) \sh (\hat{v}_{j_1} \cdots \hat{v}_{j_q}) 
= \sum_{ \alpha\in Sh(i_1\ldots i_p;i_{p+1}\ldots i_{p+q})} 
\hat{v}_{\alpha_1} \cdots \hat{v}_{\alpha_{p+q}}.
\end{equation}
By separating the tensor components and tracking the parity of each $\xi_k$ under permutation, we obtain:
\[
(\hat{v}_{i_1} \cdots \hat{v}_{i_p}) \sh (\hat{v}_{j_1} \cdots \hat{v}_{j_q}) 
= \sum_{\alpha\in Sh(i_1\ldots i_p;i_{p+1}\ldots i_{p+q})} 
(-1)^{\phi(\alpha)} \, 
\xi_{i_1} \cdots \xi_{j_q} \otimes 
v_{\alpha_1} \cdots v_{\alpha_{p+q}}.
\]

Hence, by factoring out the Grassmann algebra part, we define the \emph{super shuffle product} in $T(V)$ by:
\[
v_{i_1} \cdots v_{i_p} \shh v_{j_1} \cdots v_{j_q}
:= \sum_{\alpha\in Sh(i_1\ldots i_p;i_{p+1}\ldots i_{p+q})} 
(-1)^{\phi(\alpha)} 
v_{\alpha_1} \cdots v_{\alpha_{p+q}}.
\]

For instance,  \[
v_1 v_2 \shh v_3=v_1 v_2 v_3 +(-1)^{|v_2||v_3|}v_1 v_3 v_2+(-1)^{(|v_2|+|v_1|)|v_3|}v_3 v_1 v_2 .\]




Now, using the above notation, we adapt \cite[Proposition 1.8]{L} to the super setting. The proof follows the same arguments as in \cite[Proposition 1.8]{L}, with appropriate adjustments of the sign rule.
\ssbegin{Proposition}\label{Loday's prop} Let $v_1, v_2, \dots$ be a homogeneous basis of $V$.
The tensor algebra $(T(V),\cdot)$ with the product $\cdot,$  is a free Zinbiel superalgebra, where $\cdot$ is defined as follows: $$(v_{i_1}\cdots v_{i_p}) (v_{j_1}\cdots v_{j_q})=(v_{i_1}\cdots v_{i_p}\shh v_{j_{1}}\cdots v_{j_{q-1}})v_{j_q}.$$ Moreover,  the following set of elements 
$${B}=\cup_n\{v_{i_1}v_{i_2}\cdots v_{i_n}\,|\,n\geq 1\}$$
is a base of the free Zinbiel superalgebra $T(V).$
\end{Proposition}

For a given superset $X = X_{\bar 0} \cup X_{\bar 1},$ we consider the  free Zinbiel superalgebra $(\mathrm{Zin}(X),\circ)$ on a set of generators $X.$  For $a_1,\ldots,a_n\in \mathrm{Zin}(X)$ denote by $a_1a_2\cdots a_n$ a left-bracketed element $(\cdots(a_1\circ a_2)\cdots)\circ a_n.$ 

\ssbegin{Proposition}\label{prop: propert of sh} The super shuffle product on Zinbiel superalgebra has the following properties:

\begin{itemize}
\item[(i)] 
The super shuffle product is supercommutative and associative; namely, 
$$a\shh b=(-1)^{|a||b|}b\shh a, \quad (a\shh b)\shh c=a\shh (b\shh c), \quad \mbox{for any homogeneous} \quad a,b,c \in \mathrm{Zin}(X).$$

\item[(ii)]
$(x_{i_1}\cdots x_{i_p})\circ (x_{j_1}\cdots x_{j_q})=\begin{cases}  	
x_{i_1}\cdots x_{i_p}x_{j_1}, \quad \text{if }\quad  q=1, \\ (x_{i_1}\cdots x_{i_p}\shh x_{j_1}\cdots x_{q_{j-1}})\circ x_{j_q},  \quad \text{if }\quad q>1. \end{cases}$

\item[(iii)] Finally,

$(x_{i_1}\cdots x_{i_p})\shh (x_{j_1}\cdots x_{j_q})=$ 

$\begin{cases}  	
x_{i_1}x_{j_1}+(-1)^{|x_{i_1}||x_{j_1}|}x_{j_1}x_{i_1},  \quad \text{if }\quad p=q=1, \\[2mm]
\displaystyle (x_{i_1}\shh x_{j_1}\cdots x_{j_{q-1}})\circ x_{j_q}+(-1)^{|x_{i_1}|(\sum_{k=1}^{q}|x_{j_k}|)}x_{j_1}\cdots x_{j_q}x_{i_1},  \quad \text{if }\quad p=1 \text{ and } q>1,\\[2mm]
\displaystyle (-1)^{|x_{i_p}|( \sum_{k=1}^{q}|x_{j_k}|)}(x_{i_1}\cdots x_{i_{p-1}} \shh x_{j_1}\cdots x_{j_q})\circ x_{i_p}+(x_{i_1}\cdots x_{i_p}\shh x_{j_1}\cdots x_{j_{q-1}})\circ x_{j_q},  \text{if } p,q>1. \end{cases}$
\end{itemize}
\end{Proposition}
\begin{proof}
We will only prove Part (iii). Using the fact that $a\shh b=a\circ b+(-1)^{|a||b|}b\circ a$, we get 

\[
\begin{array}{lcl}x_{i_1} \shh (x_{j_1}\cdots x_{j_q}) & = & x_{i_1} \circ (x_{j_1}\cdots x_{j_q})+(-1)^{|x_{i_1}|(\sum_{k=1}^{q}|x_{j_k}|)} (x_{j_1}\cdots x_{j_q})\circ x_{i_1}\\[2mm]
&\stackrel{\mathrm{Part \,(ii)}}{=}& \displaystyle (x_{i_1}\shh x_{j_1}\cdots x_{j_{q-1}})\circ x_{j_q}+(-1)^{|x_{i_1}|(\sum_{k=1}^{q}|x_{j_k}|)}x_{j_1}\cdots x_{j_q}x_{i_1}.
\end{array}
\]
On the other hand, 
\begin{flalign*}
&(x_{i_1}\cdots x_{i_p})\shh (x_{j_1}\cdots x_{j_q})\\
&=(x_{i_1}\cdots x_{i_p})\circ (x_{j_1}\cdots x_{j_q})+ (-1)^{(\sum_{k=1}^{p}|x_{i_k}|)(\sum_{k=1}^{q}|x_{j_k}|)}(x_{j_1}\cdots x_{j_q})\circ (x_{i_1}\cdots x_{i_p})\\
&\stackrel{\mathrm{Part \,(ii)}}{=} (x_{i_1}\cdots x_{i_p}\shh x_{j_1}\cdots x_{j_{q-1}})\circ x_{j_q}+ (-1)^{(\sum_{k=1}^{p}|x_{i_k}|)(\sum_{k=1}^{q}|x_{j_k}|)}(x_{j_1}\cdots x_{j_q}\shh x_{i_1}\cdots x_{i_{p-1}})\circ x_{i_{p}}\\
&=(-1)^{|x_{i_p}|( \sum_{k=1}^{q}|x_{j_k}|)}(x_{i_1}\cdots x_{i_{p-1}} \shh x_{j_1}\cdots x_{j_q})\circ x_{i_p}+(x_{i_1}\cdots x_{i_p}\shh x_{j_1}\cdots x_{j_{q-1}})\circ x_{j_q}.
\qed
\end{flalign*}
\noqed

\end{proof}
\section{ Super skew-right-commutative elements in a free Zinbiel superalgebra}
This section explores a class of elements in the free Zinbiel superalgebra $\mathrm{Zin}(X)$ called super skew-right-commutative elements.  We show that every super skew-right-commutative element is a Tortkara element, extending the results of \cite{DIM19} to the super setting. Additionally, we consider the structure of the free special Tortkara superalgebra  $ST(X)$ and note that, unlike in the classical case, there are exceptional homomorphic images of  $ST(\{x, y\})$ that are not special, highlighting a key difference between the non-super and the super setting.

Consider the linear map  $
p:\mathrm{Zin}(X)\rightarrow \mathrm{Zin}(X)$ defined on the base elements
as follows (where $x_{i_1},\cdots, x_{i_{n-1}}, x_{i_n}, y,z\in X$):
$$
\begin{array}{rcl}
p(x_i) & = & -x_i,\\[2mm]
p(x_ix_j) & = & (-1)^{|x_i||x_j|}x_jx_i,\\[2mm]
p( x_{i_1}x_{i_2}\cdots x_{i_{m}}yz)&=& (-1)^{|y||z|}x_{i_1}x_{i_2}
\cdots x_{i_m}zy, \quad m\geq1.
\end{array}
$$ This map was first  introduced in \cite{DIM19} and we adapt it to the  super setting.

 Let $a\in \mathrm{Zin}(X)$ be a homogeneous element of degree greater than 1. Then  we set
$$\overline{a}=
a-p(a).$$

 Homogeneous elements of the form  $\overline{x_{i_1}\cdots x_{i_{n-1}}x_{i_n}}$, where $x_{i_1},\cdots, x_{i_{n-1}}, x_{i_n}\in X$ and $n\geq 2,$ will be called {\it super skew-right-commutative} or shortly {\it super skew-rcom} elements of $\mathrm{Zin}(X).$

 The following two lemmas are super analogs of Lemma 3.3, Lemma 3.4, and Lemma 3.5 in \cite{DIM19}. 

\ssbegin{Lemma}\label{Skw Skw} The product of super skew-right-commutative elements in Zinbiel superalgebra can be presented as follows \begin{align*}
\overline{\strut x_{i_1}\cdots x_{i_m}}\circ\overline{\strut x_{j_1}x_{j_2}} &= 
\overline{\strut x_{i_1}\cdots x_{i_m}x_{j_1}x_{j_2}} 
- (-1)^{|x_{i_{m-1}}||x_{i_m}|} \overline{\strut x_{i_1}\cdots x_{i_{m-2}}x_{i_m}x_{i_{m-1}}x_{j_1}x_{j_2}} \\
&\quad + (-1)^{|x_{i_{m}}||x_{j_1}|} (x_{i_1}\cdots x_{i_{m-1}}\shh x_{j_1})x_{i_m}x_{j_2} \\
&\quad - (-1)^{|x_{i_{m-1}}||x_{j_1}|} (x_{i_1}\cdots x_{i_{m-2}}x_{i_{m}}\shh x_{j_1})x_{i_{m-1}}x_{j_2} \\
&\quad - (-1)^{|x_{j_{2}}|(|x_{j_{1}}|+|x_{i_m}|)} (x_{i_1}\cdots x_{i_{m-1}}\shh x_{j_2})x_{i_m}x_{j_1} \\
&+ (-1)^{|x_{j_{2}}|(|x_{j_{1}}|+|x_{i_{m-1}}|)+|x_{i_{m-1}}||x_{i_m}|}  (x_{i_1}\cdots x_{i_{m-2}}x_{i_m}\shh x_{j_2})x_{i_{m-1}}x_{j_1};
\end{align*}
and \begin{align*}
\overline{\strut x_{i_1}\cdots x_{i_m}}\circ\overline{\strut x_{j_1}\cdots x_{j_n}} = & \;
\overline{\strut(x_{i_1}\cdots x_{i_m}\shh x_{j_1}\cdots x_{j_{n-2}})x_{j_{n-1}}x_{j_n}} \\
- (-1)^{|x_{i_{m}}||x_{i_{m-1}}|} & \overline{\strut(x_{i_1}\cdots x_{i_{m-2}} x_{i_m}x_{i_{m-1}}\shh x_{j_1}\cdots x_{j_{n-2}})x_{j_{n-1}}x_{j_n}} \\
 + (-1)^{|x_{i_{m}}|\sum_{k=1}^{n-1}|x_{j_{k}}|} & (x_{i_1}\cdots x_{i_{m-1}}\shh x_{j_1}\cdots x_{j_{n-1}})x_{i_m}x_{j_n} \\
 - (-1)^{|x_{i_{m-1}}|(|x_{i_{m}}|+\sum_{k=1}^{n-1}|x_{j_{k}}|)} & (x_{i_1}\cdots x_{i_{m-2}} x_{i_m}\shh x_{j_{1}}\cdots x_{j_{n-1}})x_{i_{m-1}}x_{j_n} \\
 - (-1)^{|x_{i_{m}}|(\sum_{k=1}^{n-2}|x_{j_{k}}|+|x_{j_{n}}|)+|x_{j_{n}}||x_{j_{n-1}}|}  & (x_{i_1}\cdots x_{i_{m-1}}\shh x_{j_1}\cdots x_{j_{n-2}}x_{j_n})x_{i_m}x_{j_{n-1}} \\ 
 + (-1)^{|x_{i_{m-1}}|(|x_{i_{m}}|+\sum_{k=1}^{n-2}|x_{j_{k}}|+|x_{j_{n}}|)+|x_{j_{n}}||x_{j_{n-1}}|} &(x_{i_1}\cdots x_{i_{m-2}}x_{i_m} \shh x_{j_1}\cdots x_{j_{n-2}}x_{j_n})x_{i_{m-1}}x_{j_{n-1}}.
\end{align*}
 where $m\geq2,n\geq3.$
\end{Lemma}

\begin{proof} Follows from direct calculations and Proposition \ref{prop: propert of sh}.
\end{proof}

\ssbegin{Lemma}\label{SkewRcomprudctgenerator}
The supercommutator product of super skew-rcom elements $[\overline{\strut x_{i_1}\cdots x_{i_m}},\overline{\strut x_{j_1}\cdots x_{j_n}}]$  can be presented as follows:
\begin{enumerate}
\item for $m=2$ and $q=1$ we have
$$[\overline{\strut x_{i_{1}}x_{i_{2}}},x_{j_{1}}]=\overline{\strut x_{i_1}x_{i_2}x_{j_1}}-(-1)^{|x_{i_1}||x_{i_2}|}\overline{\strut x_{i_2}x_{i_1}x_{j_1}}-(-1)^{(|x_{i_1}|+|x_{i_2}|)|x_{j_1}|}\overline{\strut x_{j_1}x_{i_1}x_{i_2}},$$

\item for $m>2$ and $q=1$ we have

\begin{align*}
[\overline{\strut x_{i_{1}}\cdots x_{i_{m-1}}x_{i_{m}}},x_{j_{1}}] = & \;
\overline{\strut x_{i_{1}}\cdots x_{i_{m}} x_{j_1}} 
- (-1)^{|x_{i_m}||x_{i_{m-1}}|} 
\overline{\strut x_{i_{1}}\cdots x_{i_{m-2}}x_{i_{m}}x_{i_{m-1}} x_{j_1}} \\
& \quad - (-1)^{|x_{j_1}|(\sum_{k=1}^{m}|x_{i_k}|)} 
\overline{\strut (x_{j_1}\shh x_{i_{1}}\cdots x_{i_{m-2}})x_{i_{m-1}}x_{i_m}},
\end{align*}

\item  for $m=2$ and $q=2$ we have  \begin{align*}
    [\overline{\strut x_{i_1}x_{i_2}},\overline{\strut x_{j_1}x_{j_2}}] = & \;
    \overline{\strut \overline{\strut x_{i_1}x_{i_2}}x_{j_1}x_{j_2}} 
    - (-1)^{(|x_{i_1}|+|x_{i_2}|)(|x_{j_1}|+|x_{j_2}|)}
    \overline{\strut \overline{\strut x_{j_1}x_{j_2}}x_{i_1}x_{i_2}} \\
  + (-1)^{|x_{i_2}||x_{j_1}|}  &  
    \overline{\strut (x_{i_1}\shh x_{j_1})x_{i_2}x_{j_2}} 
    - (-1)^{(|x_{i_2}|+|x_{j_1}|)|x_{j_2}|} 
    \overline{\strut (x_{i_1}\shh x_{j_2})x_{i_2}x_{j_1}} \\
    + (-1)^{|x_{i_1}|(|x_{i_2}|+|x_{j_2}|)+|x_{j_1}||x_{j_2}|} &
    \overline{\strut (x_{i_2}\shh x_{j_2})x_{i_1}x_{j_1}} 
    - (-1)^{|x_{i_1}|(|x_{i_2}|+|x_{j_1}|)} 
    \overline{\strut (x_{i_2}\shh x_{j_1})x_{i_1}x_{j_2}}.
\end{align*}
\item for the remaining cases we have
\begin{align*}
[\overline{\strut x_{i_1}\cdots x_{i_m}},\overline{\strut x_{j_1}\cdots x_{j_n}}] =    & \; \overline{\strut(x_{i_1}\cdots x_{i_m} \shh x_{j_1}\cdots x_{j_{n-2}})x_{j_{n-1}}x_{j_n}} \\
    - (-1)^{|x_{i_{m-1}}||x_{i_m}|} &
    \overline{\strut(x_{i_1}\cdots x_{i_{m-2}} x_{i_m}x_{i_{m-1}}\shh x_{j_1}\cdots x_{j_{n-2}})x_{j_{n-1}}x_{j_n}} \\
   - (-1)^{(|x_{i_{m-1}}|+|x_{i_m}|)\sum_{k=1}^{n}|x_{j_k}|}
     &\overline{\strut(x_{i_1}\cdots x_{i_{m-2}} \shh x_{j_1}\cdots x_{j_n} )x_{i_{m-1}}x_{i_m}} \\
    + (-1)^{(|x_{i_{m-1}}|+|x_{i_m}|)\sum_{k=1}^{n}|x_{j_k}|+|x_{j_{n-1}}||x_{j_n}|}&
    \overline{\strut(x_{i_1}\cdots x_{i_{m-2}}\shh x_{j_1}\cdots x_{j_{n-2}} x_{j_n}x_{j_{n-1}})x_{i_{m-1}}x_{i_m}} \\
   + (-1)^{|x_{i_{m-1}}|\sum_{k=1}^{n-1}|x_{j_k}|} &
    \overline{\strut(x_{i_1}\cdots x_{i_{m-1}}\shh x_{j_1}\cdots x_{j_{n-1}})x_{i_m}x_{j_n}} \\
   - (-1)^{|x_{i_{m-1}}|\sum_{k=1}^{n-1}|x_{j_k}|+|x_{i_{m-1}}||x_{i_{m}}|} &
    \overline{\strut(x_{i_1}\cdots x_{i_{m-2}} x_{i_m}\shh x_{j_1}\cdots x_{j_{n-1}})x_{i_{m-1}}x_{j_n}} \\
    - (-1)^{|x_{i_{m}}|\sum_{k=1}^{n-2}|x_{j_k}|+(|x_{i_{m}}|+|x_{j_{n-1}}|)|x_{j_{n}}|}&    \overline{\strut(x_{i_1}\cdots x_{i_{m-1}}\shh x_{j_1}\cdots x_{j_{n-2}}x_{j_n})x_{i_m}x_{j_{n-1}}} \\
    + (-1)^{|x_{i_{m}}|\sum_{k=1}^{n-2}|x_{j_k}|+(|x_{i_{m}}|+|x_{j_{n-1}}|)|x_{j_{n}}|+|x_{i_{m}}||x_{i_{m-1}}|}& \overline{\strut(x_{i_1}\cdots x_{i_{m-2}}x_{i_m} \shh x_{j_1}\cdots x_{j_{n-2}}x_{j_n})x_{i_{m-1}}x_{j_{n-1}}}.
\end{align*}
\end{enumerate}

\end{Lemma}	

\begin{proof} Follows from Proposition \ref{Loday's prop}, Proposition \ref{prop: propert of sh} and Lemma \ref{Skw Skw}.
\end{proof}

For example, for all homogeneous  $a,b,c,d\in \mathrm{Zin}(X)$ we have
\begin{equation}\label{id: L22}\begin{aligned}
[[a,b],[c,d]]=&\overline{(ab)c)d} + (-1)^{|b||c|} \overline{(ac)b)d} - (-1)^{(|b| + |c|) |d|} \overline{(ad)b)c} - (-1)^{|a||b|} \overline{(ba)c)d} \\
& - (-1)^{|a| (|b| + |c|)} \overline{(bc)a)d}  + (-1)^{|a| |b| + (|a| + |c|) |d|} \overline{(bd)a)c} \\
& + (-1)^{(|a| + |b|) |c|} \overline{(ca)b)d} - (-1)^{|a||b| + (|a| + |b|) |c|} \overline{(cb)a)d} \\
& - (-1)^{(|a| + |b|) (|c| +|d|)} \overline{(cd)a)b} - (-1)^{(|a| + |b| + |c|) |d|} (da)b)c  \\& + (-1)^{|a||b| + (|a| + |b| + |c|) |d|} \overline{(db)a)c}+
 (-1)^{(|a| + |b|) |c| + (|a| + |b| + |c|) |d|} \overline{(dc)a)b}.
\end{aligned}\end{equation}

\begin{equation}\label{id: L}
\begin{aligned}
[[[a,b],c],d]= &\overline{((ab)c)d} - (-1)^{|b||c|} \overline{((ac)b)d}  
- (-1)^{(|b| + |c|) |d|} \overline{((ad)b)c} - (-1)^{|a| |b|} \overline{((ba)c)d}  \\ &  + (-1)^{|a| (|b| + |c|)} \overline{((bc)a)d} 
+ (-1)^{|a||b| + (|a| + |c|) |d|} \overline{((bd)a)c} \\ & - (-1)^{(|a| + |b|) |c|} \overline{((ca)b)d} + (-1)^{|a| |b| + (|a| + |b|) |c|} \overline{((cb)a)d} \\ & + (-1)^{(|a| + |b|) |c| + (|a| + |b|) |d|} \overline{((cd)a)b} - (-1)^{(|a| + |b| + |c|) |d|} \overline{((da)b)c} \\ &+ (-1)^{|a| |b| + (|a| + |b| + |c|) |d|} \overline{((db)a)c} + (-1)^{(|a| + |b|) |c| + (|a| + |b| + |c|) |d|} \overline{((dc)a)b}.
\end{aligned}
\end{equation}

The following proposition is a superization of a result in \cite[Lemma 6.1]{D2}.  
\ssbegin{Proposition} Let $A=A_{\bar 0}\oplus A_{\bar 1}$ be a Zinbiel superalgebras. We have:

\textup{(}i\textup{)} The superalgebra $(A,\{\cdot , \cdot \}  )$ is supercommutative and associative.

\textup{(}ii\textup{)} The superalgebra $(A,[\cdot, \cdot])$ satisfies the Tortkara identity (\ref{id: sup-tortkara}). 
\end{Proposition}
\begin{proof} The proof is straightforward, using  \eqref{id: superZinb}, \eqref{id: L22} and \eqref{id: L}.\end{proof}

We define an element $u \in \mathrm{Zin}(X)$  as a \textit{Tortkara element} if it can be expressed by elements of $X$ in terms of super commutator product defined by  $[a, b] = ab - (-1)^{|a||b|}ba.$ Note that these elements are called  Lie elements in \cite{DIM19}. Since a Zinbiel (super)algebra under (super)commutator product is not a Lie (super)algebra,   we will refer to this element as a Tortkara element.

Let $ST(X)$ denote the free special Tortkara superalgebra generated by  $X$  under the supercommutator. This means that $ST(X)$ is a supersubalgebra of $\mathrm{Zin}(X)^{(-)} = (\mathrm{Zin}(X), [\cdot,\cdot])$  that is generated by $X.$ The degree of a word $u\in \mathrm{Zin}(X)$ will be the number of elements of the set $X$, counting repetitions, which appear in it.

Now we prove that any super skew-right-commutative element in a free superalgebra $\mathrm{Zin}(X)$ is Tortkara. The following theorem is a super analogue of Theorem 2.1 in \cite{DIM19}.

\ssbegin{Theorem}\label{Criterion}
 Let $f\in\mathrm{Zin}(X)$ be a homogeneous element of degree greater than 1. Then $f\in ST(X)$ if and only if $p(f)=-f.$  
\end{Theorem}
\begin{proof} Recall that  $ST(X)$ is generated by the supercommutator products on $X$ and $[x,y]=\overline{xy}$ for any $x,y\in X.$  By Lemma  \ref{SkewRcomprudctgenerator} we have every $f\in ST(X)$ satisfies $p(f)=-f.$

Conversely,  let us show that if $f\in \mathrm{Zin}(X)$ with $p(f)=-f$ then $f\in ST(X).$

 Write $a \equiv b$ if $a-b\in ST(X).$ If $p(f)=-f,$ then $f$ can be written as a linear combination of super skew-right-commutative elements. To establish that $f\in ST(X)$, it is sufficient to prove that \begin{equation}\label{f4}\overline{\strut x_{i_1}\cdots x_{i_n}}\equiv 0.\end{equation} We prove this by induction on $n.$ For $n=2,$ the statement is trivial. For $n=3,$ we have $$\overline{\strut x_{i_1}x_{i_2}x_{i_3}}=\frac{1}{2}[[x_{i_1},x_{i_2}],x_{i_3}]-(-1)^{|x_{i_2}||x_{i_3}|}\frac{1}{2}[[x_{i_1},x_{i_3}],x_{i_2}].$$ Now, assume that \eqref{f4} holds for all elements of degree less than $n$. That is,  $$\overline{\strut z x_{i_{k+1}}\cdots x_{i_n}}\equiv 0$$ for any 
 $z\in ST(X)$  whose degree is no more than $k.$
Setting   $z:= \overline{\strut x_{i_1}\cdots x_{i_k}}$ we obtain
$$\overline{\strut\overline{\strut x_{i_1}\cdots x_{i_k}}x_{i_{k+1}}\cdots x_{i_n}}\equiv 0$$ for $1<k<n-1.$ Hence
$$\overline{\strut x_{i_1}\cdots x_{i_{k-1}}x_{i_k}x_{i_{k+1}}\cdots x_{i_n}}\equiv(-1)^{|x_{i_{k-1}}||x_{i_k}|}\overline{\strut x_{i_1}\cdots x_{i_k}x_{i_{k-1}}x_{i_{k+1}}\cdots x_{i_n}}.$$

Since the symmetric group $S_{n-2}$ is generated by transpositions $(12),(23),\ldots,(n-3\, n-2)$,
for any $\sigma\in {S_{n-2}}$ we have 
\begin{equation}\label{f5}\overline{\strut x_{i_1}\cdots x_{i_n}}\equiv(-1)^{\phi(\sigma)}\overline{x_{\sigma(i_1)}\cdots x_{\sigma(i_{n-2})}x_{i_{n-1}}x_{i_n}},\end{equation}
where $\phi(\sigma)$ defined as in \eqref{func: phi}.

By Lemma \ref{SkewRcomprudctgenerator} we have

\begin{align*}
[\overline{\strut x_{i_{1}}\cdots x_{i_{n-2}}x_{i_{n-1}}},x_{i_{n}}]\equiv & \;
\overline{\strut x_{i_{1}}\cdots x_{i_{n-3}}x_{i_{n-2}}x_{i_{n-1}} x_{i_{n}}}-(-1)^{|x_{i_{n-2}}||x_{i_{n-1}}|}\overline{\strut x_{i_{1}}\cdots x_{i_{n-3}}x_{i_{n-1}}x_{i_{n-2}} x_{i_{n}}}\\
& - (-1)^{|x_{i_{n}}|(\sum_{k=1}^{n-1}|x_{i_k}|)} 
\overline{\strut (x_{i_n}\shh x_{i_{1}}\cdots x_{i_{n-3}})x_{i_{n-2}}x_{i_{n-1}}}\\
\equiv & \; \overline{\strut x_{i_{1}}\cdots x_{i_{n-3}}x_{i_{n-2}}x_{i_{n-1}} x_{i_{n}}}-(-1)^{|x_{i_{n-2}}||x_{i_{n-1}}|}\overline{\strut x_{i_{1}}\cdots x_{i_{n-3}}x_{i_{n-1}}x_{i_{n-2}} x_{i_{n}}}\\
& - (-1)^{|x_{i_{n}}|(\sum_{k=1}^{n-1}|x_{i_k}|)}  \sum_{\alpha\in Sh(i_{n};i_1\ldots i_{n-3})} (-1)^{\phi(\alpha)}\overline{\strut x_{i_{\alpha_1}}\cdots x_{i_{\alpha_{n-3}}}x_{i_{\alpha_{n-2}}}x_{i_{n-2}}x_{i_{n-1}}}.
\end{align*}
By \eqref{f5}, we obtain

$$[\overline{\strut x_{i_{1}}\cdots x_{i_{n-2}}x_{i_{n-1}}},x_{i_{n}}]\equiv\overline{\strut x_{i_{1}}\cdots x_{i_{n-3}}x_{i_{n-2}}x_{i_{n-1}} x_{i_{n}}}-(-1)^{|x_{i_{n-2}}||x_{i_{n-1}}|}\overline{\strut x_{i_{1}}\cdots x_{i_{n-3}}x_{i_{n-1}}x_{i_{n-2}} x_{i_{n}}}$$
$$-(-1)^{|x_{i_{n}}|(|x_{i_{n-2}}|+|x_{i_{n-1}}|)}(n-2) \overline{\strut x_{i_{1}}\cdots x_{i_{n-3}}x_{i_{n}}x_{i_{n-2}}x_{i_{n-1}}}\equiv0.
$$
Similarly, we have 
\begin{align*}[\overline{\strut x_{i_{1}}\cdots x_{i_{n-2}}x_{i_{n}}},x_{i_{n-1}}]\equiv &\; \overline{\strut x_{i_{1}}\cdots x_{i_{n-3}}x_{i_{n-2}}x_{i_{n}} x_{i_{n-1}}}-(-1)^{|x_{i_{n}}||x_{i_{n-2}}|}\overline{\strut x_{i_{1}}\cdots x_{i_{n-3}}x_{i_{n}}x_{i_{n-2}} x_{i_{n-1}}} \\
&-(-1)^{|x_{i_{n-1}}|(|x_{i_{n-2}}|+|x_{i_{n}}|)}(n-2) \overline{\strut x_{i_{1}}\cdots x_{i_{n-3}}x_{i_{n-1}}x_{i_{n-2}}x_{i_{n}}}\\
= & -(-1)^{|x_{i_{n-1}}|| x_{i_{n}}|}  \overline{\strut x_{i_{1}}\cdots x_{i_{n-2}}x_{i_{n-1}} x_{i_{n}}}-(-1)^{|x_{i_{n}}||x_{i_{n-2}}|}\overline{\strut x_{i_{1}}\cdots  x_{i_{n-3}}x_{i_{n}}x_{i_{n-2}} x_{i_{n-1}}}\\
& -(-1)^{|x_{i_{n-1}}|(|x_{i_{n-2}}|+|x_{i_{n}}|)}(n-2) \overline{\strut x_{i_{1}}\cdots x_{i_{n-3}} x_{i_{n-1}}x_{i_{n-2}}x_{i_{n}}}\equiv0.\end{align*}
Moreover, we have  
$$[\overline{\strut x_{i_{1}}\cdots x_{i_{n-3}}x_{i_{n-2}}x_{i_{n-1}}},x_{i_{n}}]+(-1)^{|x_{i_{n-1}}|| x_{i_{n}}|}
[\overline{\strut x_{i_{1}}\cdots x_{i_{n-3}}x_{i_{n-2}}x_{i_{n}}},x_{i_{n-1}}]=$$
$$-(-1)^{|x_{i_{n}}|(|x_{i_{n-2}}|+|x_{i_{n-1}}|)}(n-1) \overline{\strut x_{i_{1}}\cdots x_{i_{n-3}}x_{i_{n}}x_{i_{n-2}}x_{i_{n-1}}}$$
$$-(-1)^{|x_{i_{n-1}}||x_{i_{n-2}}|}(n-1)\overline{\strut x_{i_{1}}\cdots x_{i_{n-3}}x_{i_{n-1}}x_{i_{n-2}}x_{i_{n}}} \equiv0.$$
Thus, 

$$\overline{\strut x_{i_{1}}\cdots x_{i_{n-3}}x_{i_{n-1}}x_{i_{n-2}}x_{i_{n}}}\equiv-(-1)^{|x_{i_{n}}|(|x_{i_{n-2}}|+|x_{i_{n-1}}|)+|x_{i_{n-1}}||x_{i_{n-2}}|} \overline{\strut x_{i_{1}}\cdots x_{i_{n-3}}x_{i_{n}}x_{i_{n-2}}x_{i_{n-1}}}.$$

In other words, 
\begin{equation}\label{f6}\overline{\strut x_{i_{1}}\cdots x_{i_{n-3}}x_{i_{n-2}}x_{i_{n-1}}x_{i_{n}}}\equiv-(-1)^{|x_{i_{n}}|(|x_{i_{n-2}}|+|x_{i_{n-1}}|)+|x_{i_{n-1}}||x_{i_{n-2}}|}  \overline{\strut x_{i_{1}}\cdots x_{i_{n-3}}x_{i_{n}}x_{i_{n-1}}x_{i_{n-2}}}.\end{equation}

Set $u=x_{i_1}\cdots x_{i_{n-4}}.$ By  (\ref{f5}) and (\ref{f6}) we have  
\begin{align*} 
\overline{\strut ux_{i_{n-3}}x_{i_{n-2}}x_{i_{n-1}}x_{i_{n}}}& \equiv  -(-1)^{|x_{i_{n}}|(|x_{i_{n-2}}|+|x_{i_{n-1}}|)+|x_{i_{n-1}}||x_{i_{n-2}}|} \overline{\strut ux_{i_{n-3}}x_{i_{n}}x_{i_{n-1}}x_{i_{n-2}}}\\
& \equiv (-1)^{|x_{i_{n}}|(|x_{i_{n-3}}|+|x_{i_{n-2}}|+|x_{i_{n-1}}|)}  \overline{\strut ux_{i_{n}}x_{i_{n-3}}x_{i_{n-2}}x_{i_{n-1}}}\\
&\equiv (-1)^{|x_{i_{n}}|(|x_{i_{n-3}}|+|x_{i_{n-2}}|)+|x_{i_{n-1}}|(|x_{i_{n-3}}|+|x_{i_{n-2}}|)} \overline{\strut u x_{i_{n-1}}x_{i_{n}}x_{i_{n-3}}x_{i_{n-2}}}\\
& \equiv(-1)^{|x_{i_{n}}|(|x_{i_{n-3}}|+|x_{i_{n-2}}|)+|x_{i_{n-1}}|(|x_{i_{n-3}}|+|x_{i_{n-2}}|)+|x_{i_{n}}||x_{i_{n-1}}|} \overline{\strut ux_{i_{n}}x_{i_{n-1}}x_{i_{n-3}}x_{i_{n-2}}}\\
& \equiv (-1)^{|x_{i_{n}}||x_{i_{n-3}}|+|x_{i_{n-1}}||x_{i_{n-3}}|+|x_{i_{n}}||x_{i_{n-1}}|+|x_{i_{n-3}}||x_{i_{n-2}}|} \overline{\strut u x_{i_{n-2}}x_{i_{n}}x_{i_{n-1}}x_{i_{n-3}}}\\
&\equiv  (-1)^{|x_{i_{n}}||x_{i_{n-1}}|} \overline{\strut ux_{i_{n-3}}x_{i_{n-2}}x_{i_{n}}x_{i_{n-1}}}\\
&  \equiv - \overline{\strut u x_{i_{n-3}}x_{i_{n-2}}x_{i_{n-1}}x_{i_{n}}}.
\end{align*}

Hence $$\overline{\strut ux_{i_{n-3}}x_{i_{n-2}}x_{i_{n-1}}x_{i_{n}}}\equiv0,$$ and this completes the proof.
\end{proof} 
\begin{proof}[Another proof suggested by the referee]
Recall that $ST(X)$ is the free special Tortkara superalgebra generated by $X$  under the supercommutator (i.e.  $ST(X)$ is a supersubalgebra of $\mathrm{Zin}(X)^{(-)} = (\mathrm{Zin}(X), [\cdot,\cdot])$  that is generated by $X.$) 

We use the technique of the Grassmann envelope. 
Let $\hat{x}_j $ with the property that $\hat{x}_j= \xi_j \otimes x_j$, if  $|x_j| = \bar 1$ and $\hat{x}_j= 1 \otimes x_j$ if $|x_j|=\bar 0$ (here $\xi_j$'s are odd Grassmann variables). Denote by $\hat{X} = \{\hat{x}_j\}$ the set of such elements in the Grassmann envelope $G(\mathrm{Zin}(X))$.

Let $\Phi(x_1, \ldots, x_n)\in ST(X).$ We want to show that $\Phi(x_1, \ldots, x_n)$ can be expressed as a linear combination of monomials of the form $\overline{x_{i_1} \cdots x_{i_n}},$ where $x_{i_1}, \ldots, x_{i_n}\in X.$

Let $\hat{\Phi}(\hat{x}_1, \ldots, \hat{x}_n)$ be the corresponding commutator expression in the  Zinbiel algebra $\mathrm{Zin}(\hat{X})$. Since $\hat{x}_i$ are elements of an ordinary Zinbiel algebra, and since it is known that in the algebra case every commutator expression is a linear combination of skew-right-commutative elements (cf.  \cite[Lemma 3.6]{DIM19}), we have
\[
\hat{\Phi}(\hat{x}_1, \ldots, \hat{x}_n) = \sum_{\alpha\in S_n} \lambda_\alpha \cdot \overline{\hat{x}_{\alpha(1)} \cdots \hat{x}_{\alpha(n)}},
\]
for some scalars $\lambda_\alpha$ and permutations $\alpha\in S_n$.
Then, in $G(\mathrm{Zin}(X))$ we have 
$$\xi_1 \cdots \xi_n \otimes \Phi(x_1, \dots, x_n) =\sum_{\alpha\in S_n} (-1)^{\phi(\alpha)}\lambda_\alpha \cdot \xi_1\cdots \xi_n \otimes \overline{x_{\alpha(1)} \cdots x_{\alpha(n)}}$$

Hence, it follows that in $\mathrm{Zin}(X)$:
\[
\Phi(x_1, \ldots, x_n) = \sum_{\alpha\in S_n} (-1)^{\phi(\alpha)}\lambda_\alpha  \cdot \overline{x_{\alpha(1)} \cdots x_{\alpha(n)}}.
\]

Therefore, the original supercommutator expression $\Phi(x_1, \ldots, x_n)$ lies in the linear span of the super skew-right-commutative monomials $\overline{x_{i_1} \cdots x_{i_n}}$. Thus, every element $f \in ST(X)$ satisfies $p(f) = -f$.

Conversely, suppose that
\[
f = \overline{x_1 \cdots x_n} \in \mathrm{Zin}(X),
\]
a homogeneous expression in the generators $x_1,\ldots,x_n\in X.$ 

By Proposition 4.1, the images of the multilinear basic monomials from $\mathrm{Zin}(\hat{X})$ are linearly independent in $G(\mathrm{Zin}(X))$. Define
\[
\hat{f} := \overline{\hat{x}_1 \cdots \hat{x}_n}.
\]
Since $f$ is skew-right-commutative, $\hat{f}$ is an element of $ST(\hat{X})$ by \cite[Lemma 3.7]{DIM19}, i.e., it lies in the subspace generated by the commutators in $\mathrm{Zin}(\hat{X})$. Thus,
\[
\hat{f} = \sum_{\alpha\in S_n} \mu_\alpha \widehat{\Phi}(\hat{x}_{\alpha(1)}, \ldots, \hat{x}_{\alpha(n)}),
\]
where $\widehat{\Phi}$ is a multilinear commutator expression in the $\hat{x}_i$ and some scalars $\mu_\alpha$ and permutations $\alpha\in S_n$.

Then, the image of $\hat{f}$ in $G(\mathrm{Zin}(X))$ is
\[
\xi_1 \cdots \xi_p \otimes f =  \xi_1 \cdots \xi_p \otimes  \sum_{\alpha\in S_n} (-1)^{\phi(\alpha)}\mu_\alpha   \Phi(x_{\alpha(1)}, \dots, x_{\alpha(n)}),
\]
where $\Phi$ is the same supercommutator expression (a literal copy of $\widehat{\Phi}$), but now written in terms of the original generators $x_i$.

Since $\xi_1 \cdots \xi_p$ is a nonzero monomial in the Grassmann algebra and the tensor is zero if and only if the second tensor factor is zero, it follows that:
\[
f = \sum_{\alpha\in S_n} (-1)^{\phi(\alpha)}\mu_\alpha \Phi(x_1, \ldots, x_n) \in ST(X).
\]
Hence, $f$ belongs to the subspace generated by supercommutators, i.e., $f \in ST(X)$.
\end{proof}

It follows from Theorem \ref{Criterion} that for any $\omega\in \mathrm{Zin}(X)$ and $x\in X$ with parity of 1 we have $$\omega xx\in ST(X).$$

In \cite{DIM19}, it was proved that any Tortkara algebra with two generators is special.
We demonstrate that this is not the case in the super case. In fact,  there exists an exceptional homomorphic image of special Tortkara superalgebra $ST(\{x,y\})$ on two generators $x$ and $y,$ where $|x|=\bar 1$ and $|y|=\bar 0.$  To construct a counterexample, we use a general method that includes non-special homomorphic images of special Jordan algebras, illustrated by P.M. Cohn in \cite[Theorem 2.2]{Cohn}.  
 

Let $I$ be an ideal of $ST(X).$ By Cohn’s criterion,  $ST(X)/I$ is a 
special superalgebra if and only if $ I' \cap ST(X) \subseteq I$ where $I'$ is the ideal of a free Zinbiel superalgebra $\mathrm{Zin}(X)$ generated by
the set $I.$

\ssbegin{Theorem}\label{hom image} There exists a homomorphic image of $ST(\{x,y\})$ which is exceptional.
\end{Theorem}

\ssbegin{proof}
Let  $|x|= \bar 1$ and $|y|=\bar 0,$ and let $I$ be an ideal of $ST(\{x,y\})$ generated by elements
$$f_1=\overline{yyx},\; f_2=yxx.$$
Consider an element, $$\omega=-\overline{xyxy}-yyxx.$$
Then $ \omega=x f_1-yf_2$, thus $\omega\in I'.$ Moreover, by Theorem \ref{Criterion} we have  $\omega \in ST(\{x,y\}).$

On the other hand, it is easy to check that there are no $\lambda_1,\lambda_2\in \mathbb{K} $ such that 
$$\omega = \lambda_1 [f_1,x]+\lambda_2 [f_2,y].$$ 
 Therefore, $\omega\notin I.$ By Cohn's criterion  we conclude that  $ST(\{x,y\})/I$ is exceptional.\end{proof}

\section{Classification in low-dimension}\label{class}
Here we give a list of Tortkara superalgebras with dimensions up to three. Although the calculation is tedious, namely solving third or fourth degree equations, it is a direct calculation assisted by {\it Mathematica}. The classification is given up to an isomorphism and carried out over a field of characteristic zero. 

The following table lists all Tortkara superalgebras of dimension $2$ and $3$. The list only includes those superalgebras that are truly superalgebras in the case where $\dim=3$. We refer to  \cite{GKK18} for the list of Tortkara algebras.

In addition, we identify Tortkara superalgebras that are Malcev superalgebras. Following \cite{AE, S}, a superalgebra $M=M_{\bar 0}\oplus M_{\bar 1}$ is called a Malcev superalgebras if the following identities are satisfied: (for all  $a,b,c\in M_{\bar i}$)
\begin{eqnarray*}\label{id: sup-anti-com2}
    ab & = & -(-1)^{|a||b|}ba,\\
\label{id: sup-malcev}
    ((ab)c)d & = &-(-1)^{|a|(|b|+|c|+|d|)}((bc)d)a   -(-1)^{(|a|+|b|)(|c|+|d|)} ((cd)a)b\\[1mm]
    && -(-1)^{|d|(|a|+|b|+|c|)} ((da)b)c +(-1)^{|b||c|}(ac)(bd).
\end{eqnarray*}

 \subsubsection{$\sdim=(1|1)$} There are three classes of Tortkara superalgebras given by:

   \begin{multicols}{2}
        
        \begin{enumerate}
            \item $ {\mathtt T}_{1|1}^1=\left<e_1|e_2\right>$ \textup{(}Trivial\textup{)};
            
            \item $ {\mathtt  T}_{1|1}^2=\left<e_1|e_2;\;  e_1 \cdot e_2=e_2 \right>$;

        \columnbreak        
            \item $ {\mathtt  T}_{1|1}^3=\left<e_1|e_2; \;  e_2 \cdot e_2=e_1 \right>$.

        \end{enumerate}
    \end{multicols} 

    All of these superalgebras are Malcev types, and therefore they are Lie superalgebras by \cite{AE}.

    \begin{center}
	\begin{tabular}{|c|}
		\hline
		Malcev \\\hline

		$ {\mathtt T}_{1|1}^1$, $ {\mathtt T}_{1|1}^2$, $ {\mathtt T}_{1|1}^3$ \\\hline
        \end{tabular}

        \end{center}

 \subsubsection{$\sdim=(2|0)$} There two classes of Tortkara algebras given by:

   \begin{multicols}{2}
        
        \begin{enumerate}
            \item $ {\mathtt T}_{2|0}^1=\left<e_1|e_2\right>$ \textup{(}Trivial\textup{)};
            
            \item $ {\mathtt  T}_{2|0}^2=\left<e_1|e_2;\;  e_1 \cdot e_2=e_2 \right>$.

        \end{enumerate}
    \end{multicols} 

All these algebras are Malcev as well as Lie. 
    \subsubsection{$\sdim=(2|1)$} There are six classes of Tortkara superalgebras given in the basis $e_1, e_2|e_3$ by:
    
   \begin{multicols}{2}
        
        \begin{enumerate}
            \item $ {\mathtt T}_{2|1}^1: $ \textup{(}Trivial\textup{)};

           \item $ {\mathtt  T}_{2|1}^2: \begin{array}{lcllcl} e_3 \cdot e_3 & = & e_1
            \end{array}$

              \item $ {\mathtt  T}_{2|1}^3: \left \{\begin{array}{lcllcl} e_1 \cdot e_2 & = & - e_1; &  e_2\cdot e_3 & = &e_3;\\
            e_3\cdot e_3 & = & e_1
            \end{array} \right.$




              \item ${\bf {\mathtt  T}_{2|1}^7}: \begin{array}{lcllcl} e_1 \cdot e_2 & = & e_1
            \end{array}$

             \item $ {\mathtt  T}_{2|1}^8: \begin{array}{lcllcl} e_1 \cdot e_2 & = & e_1; &  e_1\cdot e_3 & = & e_3
            \end{array}$

             \item $ {\mathtt  T}_{2|1}^9(\gamma): \begin{array}{lcllcl} e_1 \cdot e_2 & = & \gamma e_2; &  e_1\cdot e_3 & = & e_3 
            \end{array}$






        \end{enumerate}
    \end{multicols} 
    We have $ {\mathtt  T}_{2|1}^9(\gamma)\simeq  {\mathtt  T}_{2|1}^9(\delta)$ if and only if $\gamma=\delta$.
    
Among the list above, we indicate which Tortkara superalgebras are Malcev or not.\\

\begin{center}
	\begin{tabular}{|c|}
		\hline
		Malcev \\\hline

		$ {\mathtt  T}_{2|1}^1$, $ {\mathtt  T}_{2|1}^2$, $ {\mathtt  T}_{2|1}^7$, $ {\mathtt  T}_{2|1}^9(\alpha)$ \\\hline
        \end{tabular}


\end{center}
The following superalgebras are not Lie: $T^3_{2|1},T^8_{2|1}$. 
     \subsubsection{$\sdim=(1|2)$} There are six classes of Tortkara superalgebras given in the basis $e_1, | e_2, e_3$ by:

      \begin{multicols}{2}
        
        \begin{enumerate}
            \item $ {\mathtt  T}_{1|2}^1: $ \textup{(}Trivial\textup{)};

         \item $ {\mathtt  T}_{1|2}^2: \begin{array}{lcllcl} e_1 \cdot e_3 & = &  e_2; &  e_3\cdot e_3 & = &e_1
            \end{array} $
            
                    \item $ {\mathtt  T}_{1|2}^3(\alpha): \left \{\begin{array}{lcllcl} e_1 \cdot e_2 & = &  e_2 + e_3; \\
                    e_1\cdot e_3 & = &e_2 + \alpha e_3
            \end{array} \right.$

             \item $ {\mathtt  T}_{1|2}^4(\alpha): \left \{\begin{array}{lcllcl} e_2 \cdot e_3 & = &  e_1; \\
             e_3\cdot e_3 & = &e_1; \\
             e_2 \cdot e_2 & = &\alpha e_1
            \end{array} \right.$

  \item $ {\mathtt  T}_{1|2}^5: \begin{array}{lcllcl} e_1 \cdot e_3 & = &  e_2
            \end{array} $

              \item $ {\mathtt  T}_{1|2}^6: \begin{array}{lcllcl} e_1 \cdot e_2 & = &  e_2
            \end{array}$
        
              \end{enumerate}
    \end{multicols} 
    
  We have   $ {\mathtt  T}_{1|2}^3(\alpha) \simeq  {\mathtt  T}_{1|2}^3(\beta)$ if and only if $\beta=\alpha$ or $\beta=\frac{3+\alpha}{\alpha-1}.$ Additionally, $ {\mathtt  T}_{1|2}^4(\alpha) \simeq  {\mathtt  T}_{1|2}^4(0)$ or $ {\mathtt  T}_{1|2}^3(1)$.\\

\begin{center}
	\begin{tabular}{|c|}
		\hline
		Malcev  \\\hline
		$ {\mathtt  T}_{1|2}^1$, $ {\mathtt  T}_{1|2}^2$, ${\mathtt  T}_{1|2}^3(\alpha)$, $ {\mathtt  T}_{1|2}^4(\alpha)$, $ {\mathtt  T}_{1|2}^5$, $ {\mathtt  T}_{1|2}^6$ \\\hline
        \end{tabular}

\end{center}

It was shown in \cite{AE} that every Malcev Lie superalgebra of $\dim\leq 3$ is Lie. This is not the case for Tortkara superalgebras. 
\begin{Proposition}

Every $3$-dimensional Tortkara superalgebra is a Lie superalgebra, except for ${\mathtt  T}^3_{2|1},{\mathtt  T}^8_{2|1}$.
\end{Proposition}
\begin{proof}
The proof follows from a direct computation. Let us just exhibit the computation for the algebra ${\mathtt  T}^8_{2|1}$. We have 
\[
J_s(e_1,e_2,e_3)=(e_1e_2)e_3-e_1 (e_2 e_3)-(e_1e_3)e_2=e_1e_3-0-e_3e_2=e_3-0=e_3.\qed
\]
\noqed
\end{proof}
\section*{Aknowledgement}
The authors thank Askar Dzhumadil'daev for his interest in this work. The second author is grateful to  Iryna Kashuba and Vladimir Dotsenko for useful discussions. We would like to thank the referee for contributing to the improvement of  the paper.


\end{document}